\documentclass[12pt,twoside]{article}
\usepackage{graphicx}
\usepackage{amsfonts}
\usepackage{amsbsy}
\usepackage{amssymb}
\usepackage{amsmath}
\usepackage{cite}
\usepackage{pst-all}
\usepackage{pstricks}
\usepackage{graphics}
\usepackage{tikz-cd}
\usepackage[english]{babel}
\usepackage[autostyle, english = american]{csquotes}
\MakeOuterQuote{"}

\def\bpsp{\begin{pspicture}}
\def\epsp{\end{pspicture}}
\newcommand\arrow[2]{\rotatebox{#1}{\scalebox{2}{\psline[linecolor=#2]{->}(0,0)(.03,0)}}}

\newtheorem{theorem}{Theorem}[section]
\newtheorem{remark}[theorem]{Remark}
\newtheorem{example}[theorem]{Example}
\newtheorem{lemma}[theorem]{Lemma}
\newtheorem{problem}[theorem]{Problem}
\newtheorem{corollary}[theorem]{Corollary}
\newtheorem{definition}[theorem]{Definition}
\newtheorem{proposition}[theorem]{Proposition}
\newtheorem{note}{Note}
\newtheorem{case}{Case}
\newtheorem{conjecture}{Conjecture}
\newtheorem{question}{Question}

\newcommand{\bea}{\begin{eqnarray}}
\newcommand{\eea}{\end{eqnarray}}
\newcommand{\beq}{\begin{eqnarray*}}
\newcommand{\eeq}{\end{eqnarray*}}

\def\m4{\mbox{\rm ~(mod $4$)}}

\def \bd{\begin{definition}}
\def \ed{\end{definition}}
\def \bqu{\begin{question}}
\def \equ{\end{question}}
\def \bcc{\begin{conjecture}}
\def \ecc{\end{conjecture}}
\def \bt{\begin{theorem}}
\def \et{\end{theorem}}
\def \bl{\begin{lemma}}
\def \el{\end{lemma}}
\def \bc{\begin{corollary}}
\def \ec{\end{corollary}}
\def \be{\begin{equation}}
\def \ee{\end{equation}}
\def \ben{\begin{enumerate}}
\def \een{\end{enumerate}}
\def \ba{\begin{array}}
\def \ea{\end{array}}
\def \bp{\begin{proposition}}
\def \ep{\end{proposition}}
\def \bx{\begin{example}}
\def \ex{\end{example}}
\def \br{\begin{remark}}
\def \er{\end{remark}}
\def \bdsc{\begin{description}}
\def \edsc{\end{description}}
\def \bpr{\begin{problem}}
\def \epr{\end{problem}}

\def \bn{\begin{case}}
\def \en{\end{case}}
\def \bnt{\begin{note}}
\def \ent{\end{note}}
\def\1{1\!\!1}

\def\mm2{\mbox{\rm ~(mod $2$)}}
\def\m4{\mbox{\rm ~(mod $4$)}}

\def\qed{\nolinebreak\hfill\rule{.2cm}{.2cm}\par\addvspace{.5cm}}

\def\m{\mu}

\def\1{\textbf{1}}
\def\0{\textbf{0}}

\begin{document}
\title{Essential ideals represented by mod-annihilators of modules}
\author{Rameez Raja$^a$, Shariefuddin Pirzada$^b$\\
$^{a}${\em Department of Mathematics, National Institute of Technology,}\\
{\em Hazratbal, Srinagar, India}\\
$^{b}${\em Department of Mathematics, University of Kashmir,}\\
{\em Hazratbal, Srinagar, India}\\
$^a$rameeznaqash@nitsri.ac.in; $^b$pirzadasd@kashmiruniversity.ac.in}
\date{}

\date{}
\pagestyle{myheadings} \markboth{Rameez, Pirzada}{Essential ideals represented by mod-annihilators of modules} \maketitle
\vskip 5mm
\noindent{\footnotesize \bf Abstract.}
Let $R$ be a commutative ring with unity, $M$ be a unitary $R$-module and $G$ a finite abelian group (viewed as a $\mathbb{Z}$-module). The main objective of this paper is to study properties of \textit{mod-annihilators} of $M$. For $x \in M$, we study the ideals $[x : M] =\{r\in R ~|~ rM\subseteq Rx\}$ of $R$ corresponding to \textit{mod-annihilator} of $M$. We investigate that when $[x : M]$ is an essential ideal of $R$. We prove that arbitrary intersection of essential ideals represented by \textit{mod-annihilators} is an essential ideal. We observe that $[x : M]$ is injective if and only if $R$ is non-singular and the radical of $R/[x : M]$ is zero. Moreover, if essential socle of $M$ is non-zero, then we show that $[x : M]$ is the intersection of maximal ideals and $[x : M]^2 = [x : M]$. Finally, we discuss the correspondence of essential ideals of $R$ and vertices of the annihilating graphs realized by $M$ over $R$.
\vskip 3mm

\noindent{\footnotesize Keywords: Module, ring, essential ideal, annihilator, graph.}

\vskip 3mm
\noindent {\footnotesize AMS subject classification: Primary: 13C70, 05C25.}

\section{\bf Introduction}
A nonzero ideal in a commutative ring is called essential if it intersects with every other nonzero ideal nontrivially. The study of essential ideals in a ring $R$ is a classical problem. For instance, Green and Van Wyk in \cite{GV} characterized essential ideals in certain classes of commutative and non-commutative rings. The authors in \cite{A, KR} studied essential ideals in $C(X)$, where $C(X)$ denotes the set of continuous functions on $X$. They topologically characterized the scole and essential ideals. Moreover, essential ideals have been investigated in rings of measurable functions \cite{M} and $C^{*}$- algebras \cite{KP}. For more on essential ideals, see \cite{Az, HS, J, P}.

Throughout, $R$ is a commutative ring (with $1 \neq 0$) and all modules are unitary unless otherwise stated. $[N : M] = \{r\in R~ |~ rM \subseteq N\}$ denotes an ideal of $R$. The symbols $\subseteq$  and $\subset$ have usual set theoretic meaning as containment and proper containment. We will denote the ring of integers by $\mathbb{Z}$, positive integers by $\mathbb{N}$ and the ring of integers modulo $n$ by $\mathbb{Z}_n$. For basic definitions from ring and module theory we refer to \cite{CE, W}.

For a $R$-module $M$ and $x\in M$, set $[x : M] =\{r\in R ~|~ rM\subseteq Rx\}$, which clearly is an ideal of $R$ and an annihilator of the factor module $M/Rx$, whereas the annihilator of $M$ denoted by $ann(M)$ is $[0 : M]$.

Recently in \cite{SR}, the elements of a module $M$ have been classified into \textit{full-annihilators}, \textit{semi-annihilators} and \textit{star-annihilators}. We recall a definition concerning full-annihilators, semi-annihilators and star-annihilators of a module $M$.

\begin{definition}\label{def 1} An element $x\in M$ is a:

(i) full-annihilator, if either $x = 0$ or $[x : M][y : M]M = 0$, for some nonzero $y\in M$ with $[y : M] \neq R$,

(ii) semi-annihilator, if either $x = 0$ or $[x : M] \neq 0$ and $[x : M][y : M]M = 0$, for some nonzero $y\in M$ with $ 0 \neq[y : M] \neq R$,

(iii) star-annihilator, if either $x = 0$ or $ann(M) \subset  [x : M]$ and $[x : M][y : M]M = 0$, for some nonzero $y\in M$ with  $ann(M) \subset[y : M] \neq R$.
\end{definition}

We denote by $A_f(M)$, $A_s(M)$ and $A_t(M)$ respectively the sets of full-annihilators, semi-annihilators and star-annihilators for any module $M$ over $R$ and call these annihilators as \textit{mod-annihilators}. We set $\widehat{A_{f}(M)} = A_f(M)\backslash \{0\}$, $\widehat{A_{s}(M)} = A_s(M)\backslash \{0\}$ and  $\widehat{A_t(M)} = A_t(M)\backslash \{0\}$.

This paper is organized as follows. In Section 2, we study the correspondence of essential ideals in $R$ and submodules of $M$ represented by \textit{mod-annihilators}. For some finite abelian group $G$ (viewed as a $\mathbb{Z}$-module), we determine the value of $n$ such that $[x : G] = n\mathbb{Z}$, where $x \in G$. We characterize all $\mathbb{Z}$-module $M$ such that $[x : M]$ is an essential ideal of $R$. Furthermore, we discuss that when $[x : M]$ as a $R$-module is injective and prove that if essential socle of $M$ is non-zero, then $[x : M]$ is the intersection of maximal ideals and $[x : M]^2 = [x : M]$. In Section 3, we  discuss the correspondence of essential ideals of $R$ and vertices of the annihilating graphs realized by modules over commutative rings. We conclude this paper with a discussion on some problems in this area of research.

\section{\bf Essential ideals represented by mod-annihilators}

In this section, we discuss the correspondence of essential ideals in $R$ represented by elements of $\widehat{A_{f}(M)}$, and submodules of $M$ generated by elements of $\widehat{A_{f}(M)}$. We characterize essential ideals corresponding to $\mathbb{Z}$-modules. We discuss the cases of finite abelian groups where essential ideals which are represented by elements of $\widehat{A_{f}(M)}$ corresponding to submodules of $M$ are isomorphic. If $M$ is a non-simple $R$-module, then for $x\in \widehat{A_{f}(M)}$, we show that an ideal $[x : M]$ considered as an $R$-module is injective. We also study essential ideals represented by mod-annihilators over hereditary and regular rings.

By Definition \ref{def 1}, we see that there is a correspondence of ideals in $R$ represented by elements of $\widehat{A_{f}(M)}$, $\widehat{A_{s}(M)}$, and $\widehat{A_{t}(M)}$ and cyclic submodules of $M$ generated by elements of sets $\widehat{A_{f}(M)}$, $\widehat{A_{s}(M)}$, and $\widehat{A_{t}(M)}$. Furthermore, the containment $A_t(M) \subseteq A_s(M) \subseteq A_f(M)$ is clear, so our main emphasis is on the set $\widehat{A_{f}(M)}$. However, one can study these sets separately for any module $M$.

Let $\lambda = (\lambda_1, \lambda_2, \cdots, \lambda_r)$ be a partition of $n$ denoted by $\lambda \vdash n$. For any  $\mu \vdash n$, we have an abelian group of order $p^n$ and conversely every abelian group corresponds to some partition of $n$. In fact, if $H_{\mu, p} =  \mathbb{Z}/p^{{\mu}_1}\mathbb{Z} ~\oplus~ \mathbb{Z}/p^{{\mu}_2}\mathbb{Z} ~\oplus~ \cdots ~\oplus~ \mathbb{Z}/p^{{\mu}_r}\mathbb{Z}$ is a subgroup of $G_{\lambda, p} = \mathbb{Z}/p^{\lambda_1}\mathbb{Z} \oplus \mathbb{Z}/p^{\lambda_2}\mathbb{Z} \oplus \cdots \oplus \mathbb{Z}/p^{\lambda_r}\mathbb{Z}$, then $\mu_1 \leq \lambda_1, \mu_2 \leq \lambda_2, \cdots, \mu_r \leq  \lambda_r$. If these inequalities holds we write $\mu \subset \lambda$, that is a \textquotedblleft containment order\textquotedblright on partitions. For example, a $p$-group $\mathbb{Z}/p^{5}\mathbb{Z} ~\oplus~ \mathbb{Z}/p\mathbb{Z} ~\oplus~ \mathbb{Z}/p\mathbb{Z}$ is of type $\lambda = (5, 1, 1)$. The possible types for its subgroup are: $(5, 1, 1), (4, 1, 1), (3, 1, 1), (2, 1, 1), (1, 1, 1), 2(5, 1), 2(4, 1), 2(3,1), 2(2, 1), 2(1, 1), (5), (4), (3), (2), 2$

\noindent $(1)$. Note that the types $(5, 1), (4, 1), (3,1), (2, 1), (1, 1)$ are appearing twice in the sequence of partitions for a subgroup.

Let $\lambda = (1, 1, \cdots, 1) = (1^n)$. A group of type $\lambda$ is nothing but the $\mathbb{Z}/p\mathbb{Z}$-vector space $\mathbb{Z}/p\mathbb{Z} \oplus \mathbb{Z}/p\mathbb{Z} \oplus \cdots \oplus \mathbb{Z}/p\mathbb{Z}$. Its subgroups are of type $(1^r )$, where $0\leq r \leq n$. The essential ideals corresponding to subspaces of vector space $\mathbb{Z}/p\mathbb{Z}\oplus \mathbb{Z}/p\mathbb{Z}\oplus \cdots \oplus \mathbb{Z}/p\mathbb{Z}$ 
 (represented by elements of the set $A_f(\mathbb{Z}/p\mathbb{Z}\oplus \mathbb{Z}/p\mathbb{Z} \oplus \cdots \oplus \mathbb{Z}/p\mathbb{Z})$ are same. In fact, $[x : \mathbb{Z}/p\mathbb{Z}\oplus \mathbb{Z}/p\mathbb{Z}~\oplus~ \cdots ~\oplus~ \mathbb{Z}/p\mathbb{Z}] = ann(\mathbb{Z}/p\mathbb{Z}\oplus \mathbb{Z}/p\mathbb{Z}\oplus \cdots \oplus \mathbb{Z}/p\mathbb{Z}) = p\mathbb{Z}$.

More generally, for a finite abelain $p$-group of the type $\mathbb{Z}/p^{\alpha}\mathbb{Z}\oplus \mathbb{Z}/p^{\alpha}\mathbb{Z} \cdots \oplus \mathbb{Z}/p^{\alpha}\mathbb{Z}$, where $\alpha \geq 2$. The essential ideals represented by elements of the set $A_f(\mathbb{Z}/p^{\alpha}\mathbb{Z}\oplus \mathbb{Z}/p^{\alpha}\mathbb{Z}\oplus \cdots \oplus \mathbb{Z}/p^{\alpha}\mathbb{Z}) = p^{\alpha}\mathbb{Z}$.

A finite abelian group is isomorphic to the group of the form $\mathbb{Z}/p^{\alpha_1}\mathbb{Z}\oplus \mathbb{Z}/p^{\alpha_2}\mathbb{Z}\oplus \cdots \oplus \mathbb{Z}/p^{\alpha_n}\mathbb{Z}$ whereas a finitely generated abelian group with Betti number $n$ is of the from  $\mathbb{Z}/p^{\alpha_1}\mathbb{Z}\oplus \mathbb{Z}/p^{\alpha_2}\mathbb{Z}\oplus \cdots \oplus \mathbb{Z}/p^{\alpha_n}\mathbb{Z} \oplus \mathbb{Z}\oplus \cdots \oplus \mathbb{Z}$. It is very difficult to determine the exact ideals represented by mod-annihilators of  sets $A_f(\mathbb{Z}/p^{\alpha_1}\mathbb{Z}\oplus \mathbb{Z}/p^{\alpha_2}\mathbb{Z}\oplus \cdots \oplus \mathbb{Z}/p^{\alpha_n}\mathbb{Z})$ and $A_f(\mathbb{Z}/p^{\alpha_1}\mathbb{Z}\oplus \mathbb{Z}/p^{\alpha_2}\mathbb{Z}\oplus \cdots \oplus \mathbb{Z}/p^{\alpha_n}\mathbb{Z} \oplus \mathbb{Z}\oplus \cdots \oplus \mathbb{Z})$. However, it is clear from the definition of mod-annihilators that for some $x \in A_f(\mathbb{Z}/p^{\alpha_1}\mathbb{Z}\oplus \mathbb{Z}/p^{\alpha_2}\mathbb{Z}\oplus \cdots \oplus \mathbb{Z}/p^{\alpha_n}\mathbb{Z} \oplus \mathbb{Z}\oplus \cdots \oplus \mathbb{Z})$, $[x : \mathbb{Z}/p^{\alpha_1}\mathbb{Z}\oplus \mathbb{Z}/p^{\alpha_2}\mathbb{Z}\oplus \cdots \oplus \mathbb{Z}/p^{\alpha_n}\mathbb{Z} \oplus \mathbb{Z}\oplus \cdots \oplus \mathbb{Z}]$ is some ideal in $\mathbb{Z}$.

Using the description given above, we now characterize all essential ideals represented by elements of $\widehat{A_f(M)}$ and corresponding to $\mathbb{Z}$-modules.

\begin{lemma} If $M$ is any $\mathbb{Z}$-module, then $[x : M]$ is an essential ideal if and only if $[x : M]$ is non-zero for all $x\in \widehat{A_f(M)}$.
\end{lemma}
\noindent{\bf Proof.} Let $M$ be a  $\mathbb{Z}$-module. Clearly, $M$ is an abelian group in a unique way. For all $x\in \widehat{A_f(M)}$, we have $[x : M] = n\mathbb{Z}$, $n\in\mathbb{N}$. The ideal $n\mathbb{Z}$ intersects non-trivially with any ideal $m\mathbb{Z}$, $m\in\mathbb{N}$ in $\mathbb{Z}$. So, if $M$ is a non-simple $\mathbb{Z}$-module, then for every $x\in M$, it follows that $[x : M]$ is an essential ideal. Note that $M$ is simple if and only if $\widehat{A_f(G)} = \emptyset$.

If possible, suppose $[x : M] = \{0\}$, then $[x : M]$ does not intersect non-trivially with non-trivial ideals of $\mathbb{Z}$, a contradiction.\qed
 
Since it is possible to have some finitely generated $\mathbb{Z}$-modules such that the set of mod-annihilators is equal to zero only which of course by definition is not an essential ideal. Consider a $\mathbb{Z}$-module $M = \mathbb{Z}\oplus \mathbb{Z}\oplus \cdots \oplus \mathbb{Z}$, which is a direct sum of $n$ copies of $\mathbb{Z}$. It is easy to verify that $\widehat{A_f(M)} = \widehat{M}$ with $[x : M][y : M]M = 0 $ for all $x, y\in M$. The cyclic submodules generated by elements  of $\widehat{A_f(M)}$ are simply lines with integral coordinates passing through the origin in the hyperplane $\mathbb{R}\oplus \mathbb{R}\oplus \cdots \oplus \mathbb{R}$ and these lines intersect at the origin only. Thus, for each $x \in M$, it follows that $[x : M]$ is not an essential ideal in $\mathbb{Z}$. In fact $[x : M]$ is a zero-ideal in $\mathbb{Z}$.\\ 

For any $R$-module $M$ and $x\in \widehat{A_f(M)}$, it would be interesting to characterize essential ideals $[x : M]$ represented by elements of $\widehat{A_f(M)}$ such that the intersection of all essential ideals is again an essential ideal. It is easy to see that a finite intersection of essential ideals in any commutative ring is an essential ideal. But an infinite intersection of essential ideals need not to be an essential ideal, even a countable intersection of essential ideals in general is not an essential ideal, as can be seen in \cite{A}. If the cardinality of $M$ is finite over $R$, then the submodules determined by elements of $\widehat{A_f(M)}$ are finite and therefore the ideals corresponding to submodules are finite in number. Thus, we conclude that for every $x\in \widehat{A_f(M)}$, the intersection of essential ideals $[x : M]$ in $R$ is an essential ideal. For the other case, that is, if the cardinality of $M$ is infinite over $R$, we have the following result. Note that, a nonzero submodule of a module $M$ is said to be an essential submodule of $M$ if it intersects non-trivially with other nonzero submodules of $M$.

\begin{theorem}\label{thm 1} Let $M$ be a $R$-module such that every proper submodule of $M$ is cyclic over $R$. For $x\in \widehat{A_f(M)}$, if the submodule generated by $x$ intersects non-trivially with every other nonzero submodule of $M$, then $[x : M]$ is an essential ideal in $R$.
\end{theorem}
\noindent{\bf Proof.} Assume $\bigcap\limits_{0\neq x\in M} Rx \neq 0$. If $\widehat{A_f(M)} = \phi$, then $M$ is simple, a contradiction. Let $x\in \widehat{A_f(M)}$ and let $Rx$ be the submodule generated by $x$. Since $Rx$ intersects non-trivially with every other submodule, so there exist $y\in \widehat{A_f(M)}$ such that $Rx \cap Ry \neq 0$. It suffices to prove the result for $Rx \cap Ry$. Let $z\in Rx \cap Ry$ and let $[x : M]$, $[y: M]$,    $[z : M]$ be ideals of $R$ corresponding to submodules $Rx$, $Ry$ and $Rz$. Then $[z : M]\subseteq [x : M]\cap [y: M] \neq 0$, which implies $[x : M]$ intersects non-trivially with every nonzero ideal corresponding to the submodule generated by an element of $\widehat{A_f(M)}$. For any other ideal $I$ of $R$, it is clear that $IM = \{\sum\limits_{finite}am : a\in I, ~m\in M\} = Ra$ for some $a\in M$. Thus $I$ corresponds to the cyclic submodule generated by $a\in M$. It follows that $[x : M]\cap I \neq 0$, for every nonzero ideal of $R$ and we conclude that $[x : M]$ is an essential ideal for each $x\in \widehat{A_f(M)}$.\qed

The converse of Theorem \ref{thm 1} is not true in general. We can easily construct examples from $\mathbb{Z}$-modules such that an ideal corresponding to the submodule generated by some element of $\widehat{A_f(M)}$ is an essential ideal, but the intersection of all submodules determined by elements of $\widehat{A_f(M)}$ is empty. However, if every ideal $[x : M]$, where $x\in \widehat{A_f(M)}$ corresponds to an essential submodule of $M$, then we have a non-zero intersection.

\begin{corollary} Let $M$ be a $R$-module. 

(i) For $x\in \widehat{A_f(M)}$, if the cyclic submodule $Rx$ intersects with every other cyclic nonzero submodule of $M$ non-trivially, then $[x : M]$ is an essential ideal in $R$.

(ii) The intersection $\bigcap\limits_{x\in \widehat{A_f(M)}}[x : M]$ is an essential ideal in $R$ if and only if every submodule of $M$ is essentially cyclic over $R$.
\end{corollary}

In the preceding results, we proved that  "arbitrary intersection of essentials ideals is an essential ideal". We formulated this theory of essential ideals using the concept of mod-annihilators and mainly the theory involves study of cyclic submodules of $M$. It is interesting to develop a similar theory that would employ the other finitely generated submodules of $M$. So, motivated by \cite{A}, we have the following question regarding essential ideals represented by elements of $\widehat{A_f(M_N)}$, where $\widehat{A_f(M_N)}  = \{r\in R ~|~ rM\subseteq N\}$, $N$ is a finitely generated submodule of $M$.  

\begin{problem}\label{p1} Let $M$ be a $R$-module. For $x\in \widehat{A_f(M_N)}$, characterize essential ideals $[x : M]$ in $R$ such that their intersection is an essential ideal.
\end{problem}

For a $R$-module $M$, let $Z(M)$ denote the following.
\begin{center}
$Z(M) = \{m\in M : ~ann(m) ~is ~an ~essential ~ideal ~in ~R\}$.
\end{center}
If $Z(M) = M$, then $M$ is said to be singular and if $Z(M) = 0$, then $M$ is said to be non-singular. By $rad(M)$, we denote the intersection of all maximal submodules of $M$. So, $rad(R)$ is the Jacobson radical $J(R)$ of a ring $R$. The socle of an $R$-module $M$ denoted by $Soc(M)$ is the sum of simple submodules or equivalently the intersection of all essential submodules. To say that $Soc(M)$ is an essential socle is equivalent to saying that every cyclic submodule of $M$ contains a simple submodule of $M$. An essential socle of $M$ is denoted by $essoc(M)$.

\begin{lemma}\label{lem 2} Let $M$ be a $R$-module with $essoc(M)\neq 0$, $\bigcap\limits_{0\neq x\in M} Rx \neq 0$. Then for $x\in \widehat{A_f(M)}$,  $R/[x : M]$ is a singular module.
\end{lemma}
\noindent{\bf Proof.} Since $\bigcap\limits_{0\neq x\in M} Rx \neq 0$ and $essoc(M)\neq 0$, therefore, $\widehat{A_f(M)} \neq \emptyset$. Thus, $[x : M]$ is an essential ideal. Moreover, $Z\big(R/[x : M]\big) = R/[x : M]$. Therefore, $R/[x : M]$ is a singular module.\qed

A ring $R$ is said to be a \textit{regular} ring if for all $a \in R$, $a^2x = a$ for some $x\in R$.

\begin{lemma}\label{lem 3}\cite{RZ}
A commutative ring $R$ with unity is regular if and only if every simple $R$-module is injective.
\end{lemma}

Now, we consider singular simple $R$-modules (ideals) which are injective, and obtain some properties of essential ideals corresponding to submodules generated by elements of $\widehat{A_f(M)}$.

\begin{theorem}\label{thm 3} Let $M$ be a $R$-module with $essoc(M)\neq 0$ and $\bigcap\limits_{0 \neq x\in M} Rx \neq 0$. Then every singular simple $R$-module $[x : M]$, $x\in \widehat{A_f(M)}$ is injective if and only if $Z(R) = 0$ and $rad(R/[x : M]) = 0$.
\end{theorem}
\noindent{\bf Proof.} We have $essoc(M)\neq 0$ and $\bigcap\limits_{0\neq x\in M} Rx \neq 0$, so that $\widehat{A_f(M)} \neq \emptyset$. Therefore corresponding to every cyclic submodule generated by elements of $\widehat{A_f(M)}$, we have an ideal in $R$. For $x\in \widehat{A_f(M)}$, suppose all singular simple $R$-modules $[x : M]$  are injective. If for some $z\in \widehat{A_f(M)}$, $I = [z : M]\subseteq Z(R)$ is a simple $R$-module, then $Z(I) = I$. This implies that $I$ is injective and thus a direct summand of $R$. However, the set $Z(R)$ is free from nonzero idempotent elements. Therefore, $I = 0$ and so  $Z(R) = 0$. For $x\in \widehat{A_f(M)}$, clearly $A = [x : M]$ is an essential ideal of $R$. Thus, by Lemma \ref{lem 2}, $R/A$ is a singular module and so is every submodule of $R/A$. Therefore every simple submodule of $R/A$ is injective, which implies that every simple submodule is excluded by some maximal submodule. Thus we conclude that $rad(R/A) = 0$.

For the converse, we again consider the correspondence of cyclic submodules of $M$ and ideals of $R$. Let $\tilde{I}$ be a singular simple $R$-module corresponding to the submodule of $M$. In order to show that $\tilde{I}$ is injective, we must show that for every essential ideal $A$ in $R$ corresponding to the submodule determined by an element $x\in\widehat{A_f(M)}$, every $\varphi\in Hom_R(A, \tilde{I})$ has a lift $\psi\in Hom_R(R, \tilde{I})$ such that the following diagram commutes.
\begin{center}
\begin{tikzcd}
 A \arrow{r}{i} \arrow[swap]{dr}{\varphi} & R \arrow{d}{\psi} \\
     & \tilde{I}
\end{tikzcd}
\end{center}
Let $K = ker(\varphi)$. We claim that $K$ is an essential ideal of $R$. For, if $K \cap J = \{0\}$, for some nonzero ideal $J$ of $R$, then $I^{*} = J \cap A \neq 0$ and $I^{*} \cap K = \{0\}$. This implies that $I^{*}\subseteq \varphi(I^{*}) \subseteq \tilde{I}$, a contradiction, since $\tilde{I}$ is a singular simple submodule and $Z(R) = 0$. For $\mu \neq 0$, it is clear that $\varphi$ induces an isomorphism $\mu : A/K \longrightarrow \tilde{I}$. So, $A/K$ is a simple $R$-submodule of $R/K$. By our assumption, $rad(R/K) = 0$, so there is a maximal submodule $M/K$ such that $R/K = A/K \oplus M/K$. Let $g : R \longrightarrow R/K$ be a canonical map and let $p : R/K \longrightarrow A/K$ be a projection map. Then, we have $pg : R \longrightarrow A/K$. Therefore the composition $\psi = \mu pg : R \longrightarrow \tilde{I}$ is the required lift such that the above diagram commutes.\qed

Now, we discuss some interesting consequences of the preceding theorem.

\begin{theorem}\label{thm 4} Let $M$ be an $R$-module with $essoc(M)\neq 0$, $\bigcap\limits_{0\neq x\in M} Rx \neq 0$ and for $x\in \widehat{A_f(M)}$, let every singular simple $R$-module $[x : M]$ be injective. Then every ideal $[x : M]$ is an intersection of maximal ideals, $J(R)^{2} = 0$ and $[x : M]^{2} = [x : M]$.
\end{theorem}
\noindent{\bf Proof.} For any $x\in \widehat{A_f(M)}$, clearly $[x : M]$ is an essential ideal in $R$. Therefore, $J(R) \subseteq [x : M]$, since $J(R)$ is contained in every essential ideal of $R$. On the other hand, intersection of all essential ideals in $R$ is Socle of $R$, therefore $J(R) \subseteq Soc(R)$. This implies that $J(R)^{2} = 0$ and $[x : M]$ is the intersection of maximal ideals in $R$. Suppose that $[x : M]^{2} \neq [x : M]$, for an essential ideal $[x : M]$ of $R$. By Theorem \ref{thm 3}, $Z(R) = 0$ and therefore for every essential ideal $I$, we have $I \subseteq I^2$. In particular, $[x : M] \subseteq [x : M]^2$ for each $x\in \widehat{A_f(M)}$. It follows that $[x : M]^2$ is an essential ideal and is the intersection of maximal ideals in $R$. Finally, if $y\in [x : M]^2$, $y\notin [x : M]$, there is some maximal ideal $P$ of $R$ such that $[x : M] \subseteq P$, $y\notin P$. Then $R = Ry + P$, that is, $1 = ry + m$. This implies that $y = yry + ym\in P$, a contradiction. Hence we conclude that $[x : M]^{2} = [x : M]$.\qed

\begin{corollary} Let $M$ be an $R$-module, where $R$ is hereditary. For $x\in \widehat{A_f(M)}$, if $[x : M]$ is an essential ideal of $R$ and $J(R)^2 = 0$, then every singular simple $R$-module $[x : M]$ is injective.
\end{corollary}
\noindent{\bf Proof.} Let $R$ be hereditary. From \cite{CE}, the exact sequence
\begin{equation*}
 0 \longrightarrow ann(x) \longrightarrow R \longrightarrow Rx \longrightarrow 0
\end{equation*}
splits for any $x\in R$. Since $J(R)^{2} = 0$ and $R/J(R)$ is an artinian ring, therefore $J(R) \subseteq Soc(R)$. But any essential ideal of $R$  contains $Soc(R)$. So, $J(R) \subseteq [x : M]$. This implies that $R/[x : M]$ is a completely reducible $R$-module and therefore $rad(R/[x : M]) = 0$. Thus, by Theorem \ref{thm 3}, every singular simple $R$-module $[x : M]$ is injective.\qed

Next, we consider the modules over regular rings.

\begin{theorem}\label{thm 5} Let $M$ be an $R$-module such that every submodule of $M$ is cyclic over $R$ and $\bigcap\limits_{0\neq x\in M} Rx \neq 0$. The following are equivalent.

(i) $R$ is regular

(ii) $A^2 = A$ for each ideal $A$ of $R$

(iii) $[x : M]^{2} = [x : M]$ for each $x\in \widehat{A_f(M)}$
\end{theorem}
\noindent{\bf Proof.} The equivalence of (i) and (ii) is clear and certainly (ii) implies (iii). Thus, we just need to show that (iii) implies (ii). By Theorem \ref{thm 3}, $[x : M]$ is an essential ideal for each $x\in \widehat{A_f(M)}$. Suppose $[x : M]^{2} = [x : M]$. Choose $J$ to be maximal ideal of $R$ such that $A\cap J = 0$, where $A$ is some non essential ideal of $R$. Then $A + J$ is an essential ideal of $R$. Therefore again by Theorem \ref{thm 3}, $A + J$ corresponds to some submodule of $M$ and we have $A + J = [z : M]$ for some $z\in M$. So, $(A + J)^2 = A^2 + J^2 = A + J$. If $x \in A$, then $x = \sum\limits_{finite}ab ~+  \sum\limits_{finite}mn$, where $a, b\in A$ and $m, n\in J$. Therefore,

\begin{equation*}
x - \sum\limits_{finite}ab = \sum\limits_{finite}mn \in A\cap J = 0.
\end{equation*}
This implies that $x\in A^2$ and we conclude that $A = A^2$.\qed

\begin{corollary} Let $M$ be an $R$-module with $essoc(M)\neq 0$ and $\bigcap\limits_{x\in M} Rx \neq 0$. Then every singular simple $R$-module $[x : M]$, where $x\in \widehat{A_f(M)}$, is injective if and only if $R$ is regular.
\end{corollary}
\noindent{\bf Proof.} By Theorem \ref{thm 4}, if every singular simple $R$-module $[x : M]$ is injective, then for $x\in \widehat{A_f(M)}$, we have $[x : M]^{2} = [x : M]$. Therefore, by Theorem \ref{thm 5}, $R$ is regular. If $R$ is regular, then by Lemma \ref{lem 3} every simple $R$-module is injective.\qed

\section{Representation of essential ideals by vertices of annihilating graphs}

In this section, we give a brief discussion on representation of essential ideals by vertices of graphs realized by modules over commutative rings.

A simple $graph$ $\Gamma$ consists of a $vertex~set$ $V(\Gamma)$ and an $edge~set$ $E(\Gamma)$, where an edge is an unordered pair of distinct vertices of $\Gamma$. One of the areas in algebraic combinatorics introduced by Beck \cite{Bk} is to study the interplay between graph theoretical and algebraic properties of an algebraic structure. Continuing the concept of associating a graph to an algebraic structure, another combinatorial approach of studying commutative rings was given by Anderson and Livingston in \cite{AdLs}. They associated a simple graph to a commutative ring $R$ with unity called the zero-divisor graph denoted by $\Gamma(R)$ with vertex set $Z^{*}(R) = Z(R)\backslash \{0\}$, where two distinct vertices $x,~ y\in Z^{*}(R)$ are adjacent in $\Gamma(R)$ if and only if $xy = 0$. The study of graph theoretical parameters and spectral properties in zero-divisor graphs of commutative rings are explored in \cite{AdLs, ALS, PBN, PBT, PBS, SR1}. In \cite{AdLs, SR1}, authors have discussed chromatic number, clique number and metric dimensions of zero-divisor graphs associated with finite commutative rings whereas \cite{PBN, PBS} are related to eigen values and Laplacian eigen values of zero-divisor graphs associated to finite commutative rings of type $\mathbb{Z}_n$ for $n = p^{N_1}q^{N_2}$, where $p < q$ are primes and $N_1, N_2$ are positive integers. The extension of zero-divisor graphs to non-commutative rings and semigroups can be found in \cite{FL, Rd}.

The combinatorial properties of zero-divisors discovered in \cite{Bk} have also been investigated in module theory. In \cite{SR}, the authors introduced annihilating graphs realized by modules over commutative rings known as \textit{full-annihilating, semi-annihilating and star-annihilating} graphs, denoted by $ann_f(\Gamma(M))$, $ann_s(\Gamma(M))$ and $ann_t(\Gamma(M))$. The vertices of annihilating graphs are elements of sets $\widehat{A_{f}(M)}$, $\widehat{A_{s}(M)}$ and $\widehat{A_{t}(M)}$ respectively, where two vertices $x$ and $y$ are adjacent if and only if $[x : M][y : M]M = 0$. The three simple graphs: full-annihilating, semi-annihilating and star-annihilating with vertex sets: $\widehat{A_f(M)}, \widehat{A_s(M)}, \widehat{A_t(M)}$ are natural generalizations of the zero-divisor graph introduced in \cite{AdLs}. This concept was further studied in \cite{SR2}.

We call a vertex $x$, an \textit{essential vertex} in $ann_f(\Gamma(M))$ if the ideal represented by $x$ is essential in $R$. Recall that a graph $\Gamma$ is said to be a \textit{complete} if there is an edge between every pair of distinct vertices.

By Definition \ref{def 1}, we see the containment $ann_t(\Gamma(M)) \subseteq ann_s(\Gamma(M)) \subseteq ann_f(\Gamma(M))$ as induced subgraphs of the graph $ann_f(\Gamma(M))$, since $A_t(M) \subseteq A_s(M) \subseteq A_f(M)$. If $ann_f(\Gamma(M))$ is a finite graph, then by [\cite{SR}, Theorem 3.3 and Example 2.2], $|\widehat{A_{f}(M)}| = |\widehat{A_{s}(M)}|$ and annihilating graphs $ann_f(\Gamma(M))$, $ann_s(\Gamma(M))$ coincide, whereas the graph $ann_t(\Gamma(M))$ with vertex set $\widehat{A_{t}(M)}$ may be different. For a $\mathbb{Z}$-module $M = \mathbb{Z}\oplus \cdots \oplus \mathbb{Z}$, we have by Definition \ref{def 1} $[x : M][y : M]M = 0$ for all $x, y\in \widehat{A_f(M)}$. Therefore $ann_f(\Gamma(M))$ is a complete graph whereas the graph $ann_s(\Gamma(M))$ is an empty graph. Thus for finitely generated infinite modules, graphs $ann_f(\Gamma(M))$ and $ann_s(\Gamma(M))$ are different.

As discussed in Section 2, for a module $M = \mathbb{Z}\oplus \cdots \oplus \mathbb{Z}$, the ideal $[x : M]$ represented by a vertex $x\in \widehat{A_f(M)}$ of the graph $ann_f(\Gamma(M))$ is not an essential ideal. So, $x$ is not an essential vertex of the graph $ann_f(\Gamma(M))$. On the other hand, every vertex of a $\mathbb{Z}$-module $\mathbb{Z}_p \oplus \mathbb{Z}_q$ is an essential vertex of the graph $ann_f(\Gamma(\mathbb{Z}_p \oplus \mathbb{Z}_q))$, where $p$ and $q$ are any two primes.\\

Finally, Problem  \ref{p1} can be restated in the graph theoretical version as follows.
\begin{problem} Characterize all annihilating graphs realized by a module $M$ such that every vertex $x \in \widehat{A_f(M_N)}$ of an annihilating graph is an essential vertex.
\end{problem}

{\bf Conclusion:} In this paper, we formulated a new approach of recognition of essential ideals  in a commutative ring $R$. This formulation of essential ideals corresponds to mod-annihilators of a $R$-module $M$. It is interesting to characterize essential ideals such that their arbitrary intersection is an essential ideal, since it is specified in \cite{A} that an arbitrary intersection of essential ideals may not be an essential ideal. Furthermore, we obtained the results related to ideals $[x : M]$ of $R$, where $x$ is a mod-annihilator of $M$ and discussed the representation of vertices of annihilating graphs by essential ideals of $R$. Apart from the research problems which we mentioned in Sections 2 and 3, the following problems could be investigated for the future work.\\

\textbf{1.} If $G$ is a finite abelian $p$-group (viewed as a finite $\mathbb{Z}$-module) of rank at least $3$. Determine value of $n$ for the essential ideal $[x : G] = n\mathbb{Z}$, where $x \in G$.

\textbf{2.} If $G$ is any finite abelian group (viewed as a finite $\mathbb{Z}$-module). Determine value of $n$ for the essential ideal $[x : G] = n\mathbb{Z}$, where $x \in G$.

\end{document}